\documentclass[10pt]{amsart}
\usepackage{amssymb}
\usepackage{amscd}

\newcounter{TmpEnumi}

\theoremstyle{definition}
\newtheorem{thm}{Theorem}[section]
\newtheorem{lem}[thm]{Lemma}
\newtheorem{prp}[thm]{Proposition}
\newtheorem{dfn}[thm]{Definition}
\newtheorem{cor}[thm]{Corollary}

\newtheorem{rmk}[thm]{Remark}
\newtheorem{ntn}[thm]{Notation}
\newtheorem{exa}[thm]{Example}

\newtheorem{qst}[thm]{Question}

\newcommand{\beq}{\begin{equation}}
\newcommand{\eeq}{\end{equation}}
\newcommand{\beqr}{\begin{eqnarray*}}
\newcommand{\eeqr}{\end{eqnarray*}}
\newcommand{\bal}{\begin{align*}}
\newcommand{\eal}{\end{align*}}
\newcommand{\bei}{\begin{itemize}}
\newcommand{\eei}{\end{itemize}}
\newcommand{\limi}[1]{\lim_{{#1} \to \infty}}

\newcommand{\af}{\alpha}
\newcommand{\bt}{\beta}
\newcommand{\gm}{\gamma}
\newcommand{\dt}{\delta}
\newcommand{\ep}{\varepsilon}

\newcommand{\ld}{\lambda}

\newcommand{\kp}{\kappa}
\newcommand{\ph}{\varphi}
\newcommand{\ps}{\psi}
\newcommand{\rh}{\rho}
\newcommand{\om}{\omega}
\newcommand{\ta}{\tau}

\newcommand{\Z}{{\mathbb{Z}}}

\newcommand{\C}{{\mathbb{C}}}
\newcommand{\N}{{\mathbb{N}}}

\newcommand{\id}{{\mathrm{id}}}

\newcommand{\dist}{{\mathrm{dist}}}

\newcommand{\Aut}{{\mathrm{Aut}}}
\newcommand{\Ad}{{\mathrm{Ad}}}

\newcommand{\card}{{\mathrm{card}}}

\newcommand{\Ker}{{\mathrm{Ker}}}

\newcommand{\dirlim}{\varinjlim}
\newcommand{\Mi}{M_{\infty}}

\newcommand{\andeqn}{\,\,\,\,\,\, {\mbox{and}} \,\,\,\,\,\,}
\newcommand{\QED}{\rule{0.4em}{2ex}}

\newcommand{\ca}{C*-algebra}

\newcommand{\pj}{projection}

\newcommand{\hm}{homomorphism}
\newcommand{\wolog}{without loss of generality}

\newcommand{\ifo}{if and only if}

\newcommand{\mops}{mutually orthogonal \pj s}
\newcommand{\pisca}{purely infinite simple \ca}

\newcommand{\hsa}{hereditary subalgebra}

\renewcommand{\S}{\subset}

\newcommand{\I}{\infty}

\title[Crossed products by Rokhlin actions]{Crossed
  products by finite group actions with the Rokhlin property}

\date{6~Feb.\  2009}

\author{Hiroyuki Osaka}

\address{Department of Mathematical Sciences,
 Ritsumeikan University, Kusatsu, Shiga, 525-8577 Japan}

\email[]{osaka@se.ritsumei.ac.jp}

\author{N.~Christopher Phillips}

\address{Department of Mathematics, University of Oregon,
      Eugene OR 97403-1222, USA.}

\email[]{ncp@darkwing.uoregon.edu}
\subjclass[2000]{Primary 46L55; Secondary 46L35.}
\thanks{Research of the first author partially supported by
The Open Research Center Project for Private Universities:
matching fund from MEXT, 2004-2008.
Research of the second author partially supported by
NSF grants DMS 0302401 and DMS 0701076.}

\begin{document}

\begin{abstract}
We prove that a number of classes of separable unital \ca s
are closed under crossed products by finite group actions with
the Rokhlin property, including:
\begin{itemize}
\item
AI~algebras, AT~algebras, and related classes characterized
by direct limit decompositions using semiprojective
building blocks.
\item
Simple unital AH~algebras with slow dimension growth
and real rank zero.
\item
\ca s with real rank zero or stable rank one.
\item
Simple \ca s for which the order on \pj s is determined by traces.
\item
\ca s whose quotients all satisfy the
Universal Coefficient Theorem.
\item
\ca s with a unique tracial state.
\end{itemize}
Along the way, we give a systematic treatment of the derivation
of direct limit decompositions from local approximation conditions
by homomorphic images which are not necessarily injective.
\end{abstract}

\maketitle

\setcounter{section}{-1}

\section{Introduction}\label{Sec:Intro}

\indent
The purpose of this paper is to prove theorems of the following type.
Let $A$ be a separable unital \ca,
let $G$ be a finite group, and let $\af \colon G \to \Aut (A)$
be an action which has the Rokhlin property.
(See Definition~\ref{SRPDfn} below.)
Suppose that $A$ belongs to a class of
\ca s characterized by some structural property,
such as the AI~algebras.
Then $C^* (G, A, \af)$ belongs to the same class.
The classes we consider include:
\begin{itemize}
\item
\ca s with various kinds
of direct limit decompositions involving semiprojective
building blocks.
\item
Simple unital AH~algebras with slow dimension growth
and real rank zero.
\item
\ca s with real rank zero or stable rank one.
\item
Simple \ca s for which the order on \pj s is determined by traces.
\item
\ca s $A$ such that all quotients $A / I$ satisfy the
Universal Coefficient Theorem.
\item
\ca s with a unique tracial state.
\end{itemize}
(For a complete list,
see the statements of the theorems
in Sections \ref{Sec:CrPrdRP} and~\ref{Sec:CrPrdRP2}.)
In most of these classes, the \ca s need not be simple.

The essential observation
(except for the last item)
is already in the proof
of Theorem~2.2 of~\cite{PhtRp1a}:
the crossed product $C^* (G, A, \af)$ has a local
approximation property by \ca s stably isomorphic to
homomorphic images of $A.$
For some of the classes we consider, this is already enough.
For others, such as classes of direct limits
involving semiprojective building blocks,
technical difficulties arise because
we must apparently allow arbitrary homomorphic images
in the approximation property.
(Homomorphic images of semiprojective \ca s need not be
semiprojective.)
We therefore present a systematic framework for proving that
a \ca\  with a local approximation property by
homomorphic images of a suitable class of semiprojective \ca s
can in fact be written as a direct limit of algebras in the class.
Our result is an analog of Lemma~15.2.2 of~\cite{Lr},
where the local approximation property
uses injective homomorphic images.
Results of this kind for specific cases are already implicit in
the literature;
our contribution is primarily to systematize the method.

We know of only three results like those above in the literature:
for AF~algebras (Theorem~2.2 of~\cite{PhtRp1a}),
for approximately divisible \ca s,
and for $D$-absorbing \ca s for appropriate $D$
(both in Corollary~3.4 of~\cite{HW}).
(Some others, for example involving the ideal and projection
properties, will appear in~\cite{PP}.
Some of them were proved before this paper was written,
but do not use the ideas of this paper.)
It is surprising to us that Theorem~2.2 of~\cite{PhtRp1a}
had not been proved long ago,
and even in~\cite{PhtRp1a} the much wider applicability
of the argument was overlooked.

Our main applications will appear in~\cite{OP4},
where we use results of the type presented here to
show that various actions of finite groups do not have the
Rokhlin property.
(Indeed, the Rokhlin property is rather rare.)
We present the results of this paper separately
for three reasons.
First, the main ideas are somewhat different from those
in~\cite{OP4}.
Second, some of the consequences,
such as those involving the Universal Coefficient Theorem
and direct limits using infinite building blocks,
have no connection with~\cite{OP4}.
Third, some of the results have other applications.
See~\cite{PV}.

In Section~\ref{Sec:DLim}, we introduce
flexible classes of \ca s,
and prove that if a separable unital \ca\  has local
approximation by homomorphic images of \ca s in a flexible class,
then it is a direct limit of algebras in the class.
The algebras involved need not be simple.
In Section~\ref{Sec:ExFlex},
we prove that many of the standard classes of
semiprojective building blocks are flexible.
The most complicated case is
the one dimensional noncommutative CW~complexes of~\cite{ELP}.
Sections \ref{Sec:CrPrdRP} and~\ref{Sec:CrPrdRP2}
contain the proofs of the
closure results under crossed products.
They are roughly divided into membership in various classes
considered in the Elliott program (Section~\ref{Sec:CrPrdRP})
and more general properties such as stable rank
(Section~\ref{Sec:CrPrdRP2}).

We take $\N$ to be $\{ 1, 2, \ldots \}.$

\section{Direct limits and local approximation}\label{Sec:DLim}

\indent
The purpose of this section is to give a uniform description of
the local approximation characterizations
of the classes of direct limit algebras constructed using some common
families of semiprojective building blocks.
We do not require our direct limits to either be simple or have
real rank zero, and we do not require the maps in our direct
systems or local approximations to be injective.
The results are known for various specific classes of building blocks.

\begin{dfn}\label{D:FSat}
Let ${\mathcal{C}}$ be a class of separable unital C*-algebras.
Then ${\mathcal{C}}$ is {\emph{finitely saturated}}
if the following closure conditions hold:
\begin{enumerate}
\item\label{D:FSat:1}
If $A \in {\mathcal{C}}$ and $B \cong A,$ then $B \in {\mathcal{C}}.$
\item\label{D:FSat:2}
If $A_1, A_2, \ldots, A_n \in {\mathcal{C}}$ then
$\bigoplus_{k=1}^n A_k \in {\mathcal{C}}.$
\item\label{D:FSat:3}
If $A \in {\mathcal{C}}$ and $n \in \N,$
then $M_n (A) \in {\mathcal{C}}.$
\item\label{D:FSat:4}
If $A \in {\mathcal{C}}$ and $p \in A$ is a nonzero projection,
then $p A p \in {\mathcal{C}}.$
\end{enumerate}
Moreover,
the {\emph{finite saturation}} of a class ${\mathcal{C}}$ is the
smallest finitely saturated class which contains ${\mathcal{C}}.$
\end{dfn}

If one wants to deal with sets instead of classes, one should
omit Condition~(\ref{D:FSat:1}).

\begin{dfn}\label{D:Flex}
Let ${\mathcal{C}}$ be a class of separable unital C*-algebras.
We say that ${\mathcal{C}}$ {\emph{has approximate quotients}} if:
\begin{enumerate}
\item\label{D:Flex:2}
For every $A \in {\mathcal{C}}$ and every ideal $I \subset A,$
there is an
increasing sequence $I_0 \subset I_1 \subset \cdots$ of ideals in
$A$ such that ${\overline{\bigcup_{n = 0}^{\I} I_n}} = I,$ and such
that for every $n$ the C*-algebra $A / I_n$ is in the finite
saturation of ${\mathcal{C}}.$
\setcounter{TmpEnumi}{\value{enumi}}
\end{enumerate}
We say that ${\mathcal{C}}$ is {\emph{flexible}} if in addition:
\begin{enumerate}
\setcounter{enumi}{\value{TmpEnumi}}
\item\label{D:Flex:1}
For every $A \in {\mathcal{C}},$ every $n \in \N,$
and every nonzero projection
$p \in M_n (A),$ the corner $p M_n (A) p$ is finitely generated,
and is semiprojective in the
sense of Definition 14.1.3 of~\cite{Lr}.
\end{enumerate}
\end{dfn}

We would like to replace Condition~(\ref{D:Flex:1}) by the simpler
requirement that all algebras in ${\mathcal{C}}$ be semiprojective.
However, we do not know whether corners in semiprojective algebras
are necessarily semiprojective (although this is true for full
corners, by Proposition~2.27 of \cite{Bl2}).

The finite generation condition is technically important in our
proofs, but we do not know if it is really necessary.
It is unfortunately not true that a corner in a unital
finitely generated \ca\  is again finitely generated.
Let $X$ be the Hilbert cube $[0, 1]^{\N},$
and let $A = [K \otimes C (X)]^+.$
By Theorem~8 of~\cite{OZ},
the algebra $K \otimes C (X)$ is generated by a single element,
so $A$ is finitely generated.
However, $A$ has a corner isomorphic to $C (X).$
If $C (X)$ were generated by $n$ elements,
we would have an injective map $X \to \C^n.$

The collection of all finite dimensional C*-algebras is finitely
saturated and flexible.
Further examples will be given in Section~\ref{Sec:ExFlex}.

\begin{lem}\label{L:FlexSat}
Let ${\mathcal{C}}$ be a class of separable unital C*-algebras.
If ${\mathcal{C}}$ is flexible,
then its finite saturation is also flexible.
If ${\mathcal{C}}$ satisfies
either condition of Definition~\ref{D:Flex} separately,
then so does its finite saturation.
\end{lem}

\begin{proof}
We can obtain the finite saturation of ${\mathcal{C}}$
by first forming all
finite matrix algebras over all algebras in ${\mathcal{C}},$
then forming
all corners of the algebras we have so far, then forming all
finite direct sums of the resulting algebras, and finally closing
under the isomorphism relation.

It is easy to prove that these operations all preserve
Condition~(\ref{D:Flex:1}) of Definition~\ref{D:Flex},
using the
following three facts: a finite direct sum of semiprojective
C*-algebras is semiprojective (Theorem 14.2.1 of \cite{Lr}),
a finite direct sum of finitely generated C*-algebras is finitely
generated (easy),
and  a corner in a finite direct sum is a finite direct sum of
corners in the summands.

Condition~(\ref{D:Flex:2}) of Definition~\ref{D:Flex}
is obviously preserved under
finite direct sums and tensoring with $M_n.$
We therefore need only show that it 
passes to corners.

Accordingly, let $A$ be a separable unital C*-algebra satisfying
Condition~(\ref{D:Flex:2}),
let $p \in A$ be a nonzero projection, and let
$I \subset p A p$ be an ideal.
Define $J = {\overline{A I A}},$ which is an ideal in $A.$
Note that
\[
p J p = {\overline{(p A) (p I p) (A p)}}
 = {\overline{(p A p) I (p A p)}}
 \subset I.
\]
It follows that $J \cap p A p = p J p = I.$

By hypothesis, there are ideals $J_0 \subset J_1 \subset \cdots$ in
$A$ such that ${\overline{\bigcup_{n=0}^{\I} J_n}} = J$ and such that
$A / J_n$ is in the finite saturation of ${\mathcal{C}}$ for all $n.$
Define $I_n = p J_n p,$ which is an ideal in $p A p.$
Note that $I_n = J_n \cap p A p \subset J \cap p A p = I.$
Also,
\[
{\textstyle{
{\overline{\bigcup_{n=0}^{\I} I_n}}
 = {\overline{\bigcup_{n=0}^{\I} p J_n p}}
 = p \left({\overline{\bigcup_{n=0}^{\I} J_n}} \right) p
 = p J p
 = I.
}}
\]

Let $\pi \colon A \to A / J_n$ and $\kp \colon p A p \to p A p / I_n$
be the quotient maps.
We next claim that $p A p / I_n \cong \pi (p) (A / J_n ) \pi (p).$
Since $J_n \cap p A p = I_n,$ we have an obvious injective map
$\ph \colon p A p / I_n \to A / J_n.$
For surjectivity, let $x \in \pi (p) (A / J_n ) \pi (p).$
Choose $a \in A$ such that $\pi (a) = x.$
Then $a - p a p \in J_n.$
So $\ph (\kp (p a p)) = \pi (p a p) = \pi (a) = x.$
This proves the claim.

Since $p A p / I_n$ is a corner of $A / J_n,$ and since $A / J_n$ is in
the finite saturation of ${\mathcal{C}},$
it follows that $p A p / I_n$ is in the
finite saturation of ${\mathcal{C}}.$
This completes the verification that $p A p$ satisfies
Condition~(\ref{D:Flex:2}), and the proof of the lemma.
\end{proof}

\begin{dfn}\label{D:AC}
Let ${\mathcal{C}}$ be a class of separable unital C*-algebras.
A {\emph{unital approximate ${\mathcal{C}}$-algebra}}
is a C*-algebra $A$ which
is isomorphic to a direct limit $\dirlim A_n$ in which each $A_n$
is in the finite saturation of ${\mathcal{C}}$
and in which each homomorphism
$\ph_n \colon A_n \to A_{n + 1}$ is unital
(but not necessarily injective).
\end{dfn}

\begin{dfn}\label{D:LC}
Let ${\mathcal{C}}$ be a class of separable unital C*-algebras.
A {\emph{unital local ${\mathcal{C}}$-algebra}}
is a separable unital C*-algebra $A$
such that for every finite set $S \subset A$ and every $\ep > 0,$
there is a C*-algebra $B$ in the finite saturation of ${\mathcal{C}}$
and a unital homomorphism $\ph \colon B \to A$
(not necessarily injective)
such that $\dist (a, \, \ph (B)) < \ep$ for all $a \in S.$
If one can always choose $B \in {\mathcal{C}},$
rather than merely in its finite saturation,
we call $A$ a {\emph{unital strong local ${\mathcal{C}}$-algebra}}
\end{dfn}

Several of the following closure properties will be needed later.
Note that, if ${\mathcal{C}}$
has approximate quotients (Definition~\ref{D:Flex}),
then the class of unital local ${\mathcal{C}}$-algebras
is finitely saturated.

\begin{lem}\label{L:ClosULC}
Let ${\mathcal{C}}$ be a
class of separable unital C*-algebras.
Let ${\mathcal{L}}$ be the class
of unital local ${\mathcal{C}}$-algebras.
Then:
\begin{enumerate}
\item\label{L:ClosULC:1}
If $A_1, A_2, \ldots, A_n \in {\mathcal{L}}$ then
$\bigoplus_{k=1}^n A_k \in {\mathcal{L}}.$
\item\label{L:ClosULC:2}
If $A \in {\mathcal{L}}$ and $n \in \N,$
then $M_n (A) \in {\mathcal{L}}.$
\item\label{L:ClosULC:3}
If $A \in {\mathcal{L}}$
and $I$ is an ideal in~$A,$
then $A / I \in {\mathcal{L}}.$
\item\label{L:ClosULC:3a}
If $A$ is a unital strong local ${\mathcal{L}}$-algebra,
then $A \in {\mathcal{L}}.$
\setcounter{TmpEnumi}{\value{enumi}}
\end{enumerate}
Suppose, in addition, that ${\mathcal{C}}$
has approximate quotients (Definition~\ref{D:Flex}).
Then:
\begin{enumerate}
\setcounter{enumi}{\value{TmpEnumi}}
\item\label{L:ClosULC:5}
If $A \in {\mathcal{L}}$
and $p \in A$ is a nonzero \pj,
then $p A p \in {\mathcal{L}}.$
\item\label{L:ClosULC:4}
If $A$ is a unital local ${\mathcal{L}}$-algebra,
then $A \in {\mathcal{L}}.$
\end{enumerate}
\end{lem}

\begin{proof}
The first four parts are obvious.

For Part~(\ref{L:ClosULC:5}),
\wolog\  ${\mathcal{C}}$ is finitely saturated.
Let $A \in {\mathcal{L}},$
let $p \in A$ be a nonzero \pj,
let $S \S p A p$ be finite,
and let $\ep > 0.$
Set $M = \sup_{a \in S} \| a \|.$
Choose $\dt > 0$ with
\[
\dt \leq \min \left( 1, \, \frac{\ep}{2 M + 3} \right)
\]
and so small that whenever $A$ is a unital \ca,
$B \S A$ is a subalgebra,
and $p \in A$ is a \pj\  such that $\dist (p, B) < \dt,$
then there are a \pj\  $q \in B$ and a unitary $u \in A$
such that
\begin{equation}\label{Eq:uq}
u q u^* = p
\andeqn
\| u - 1 \| < \frac{\ep}{2 M + 3}.
\end{equation}
Write $S = \{ a_1, a_2, \ldots, a_n \}.$
Choose $B_0 \in {\mathcal{C}}$
and a unital \hm\  $\rh_0 \colon B_0 \to A$
such that there are $x, c_1, c_2, \ldots, c_n \in B_0$
with $\| \rh_0 (x) - p \| < \dt$
and $\| \rh_0 (c_j) - a_j \| < \dt$ for $j = 1, 2, \ldots, n.$
By the choice of $\dt,$
there are a \pj\  $q \in \rh_0 (B_0)$ and a unitary $u \in A$
satisfying~(\ref{Eq:uq}).

Since ${\mathcal{C}}$ has approximate quotients,
there is an
increasing sequence $I_0 \subset I_1 \subset \cdots$ of ideals in
$A$ such that ${\overline{\bigcup_{n = 0}^{\I} I_n}} = \ker (\rh_0),$
and such
that for every $n$ we have $A / I_n \in {\mathcal{C}}.$
By semiprojectivity of $\C$
(see Definition 14.1.3 of~\cite{Lr}),
there exists $n$ and a \pj\  $f \in B_0 / I_n$
whose image in $A$ is~$q.$
Set $B = B_0 / I_n,$ which is in ${\mathcal{C}},$
let $\pi \colon B_0 \to B$ be the quotient map,
and let $\rh \colon B \to A$ be the unital \hm\  such that
$\rh_0 = \rh \circ \pi.$
Define $\ph \colon f B f \to p A p$ by
$\ph (b) = u \rh (b) u^*.$
Note that $\ph (b)$ really is in $p A p,$
because $u \rh (f) u^* = u q u^* = p.$
Set $b_j = f \pi (c_j) f$ for $j = 1, 2, \ldots, n.$
Then
\begin{align*}
\| \ph (b_j) - a_j \|
& = \| u q \rh_0 (c_j) q u^* - a_j \|
  = \| p u \rh_0 (c_j) u^* p - a_j \|
  \leq \| u \rh_0 (c_j) u^* - a_j \|
          \\
& \leq 2 \| u - 1 \| \cdot \| \rh_0 (c_j) \| + \| \rh_0 (c_j) - a_j \|
          \\
& < 2 \left( \frac{\ep}{2 M + 3} \right) (M + \dt) + \dt
  \leq 2 \left( \frac{\ep}{2 M + 3} \right) (M + 1)
       + \frac{\ep}{2 M + 3}
  = \ep.
\end{align*}
This shows that $p A p \in {\mathcal{L}}.$

For the last part,
we now know that the hypotheses imply
that ${\mathcal{L}}$ is finitely saturated.
Apply Part~(\ref{L:ClosULC:3a}).
\end{proof}

Obviously, every unital approximate ${\mathcal{C}}$-algebra
is a unital local ${\mathcal{C}}$-algebra.
%
We will prove that if ${\mathcal{C}}$ is flexible,
then every local ${\mathcal{C}}$-algebra
is an approximate ${\mathcal{C}}$-algebra.
This is related to, but not quite the same as, the results of
Section~15.2 of~\cite{Lr}.
The main difference is the issue of injectivity of the maps involved.

In a number of specific cases,
this result is known, or is at least implicit in the literature.
For example, if ${\mathcal{C}} = \{ \C, \, C ([0,1]), \, C (S^1) \}$
(see Example~\ref{E:AT} below),
it is essentially Theorem~4.1.5 of~\cite{Lnbk}.
There is a slight difference because the definition of a circle
algebra used there (Definition 2.3.1 of~\cite{Lnbk}) allows
algebras such as the continuous functions on the Cantor set.
Lemma 4.1.2 of~\cite{Lnbk} implies that every unital AT~algebra as
defined there is a unital local ${\mathcal{C}}$-algebra
in our sense.

\begin{prp}\label{P:LCFlex}
Let ${\mathcal{C}}$ be a flexible class
of separable unital C*-algebras.
Then every unital local ${\mathcal{C}}$-algebra (Definition~\ref{D:LC})
is
a unital approximate ${\mathcal{C}}$-algebra (Definition~\ref{D:AC}).
\end{prp}

\begin{proof}
By Lemma~\ref{L:FlexSat}, without loss of generality we may assume
${\mathcal{C}}$ is finitely saturated.

Let $\{ x_0, x_1, \ldots \} \subset A$ be a countable dense set.

For the purposes of this proof, if $A$ is a unital C*-algebra,
$G \subset A,$ and $n \in \N \cup \{ 0 \},$ we define $S_n (G)$
to be the set of all noncommuting polynomials in the elements
$x, x^*$ for $x \in G,$ with coefficients in
$\Z [2^{-n}] + i \Z [2^{-n}]$ and of absolute value at most $2^n.$
Then $S_n (G) \subset S_{n + 1} (G).$
If $G$ is finite, then so is $S_n (G).$
Also, if $G$ generates $A,$ then $\bigcup_{n = 0}^{\I} S_n (G)$
is dense in $A.$

We also need finite presentations of our algebras, in the sense of
Definition~13.2.2 of~\cite{Lr}.
We may always require finite generating sets to be bounded by $1.$
Thus, following Definition~13.2.1 of~\cite{Lr}, for a finite set $G$ of
generators we consider the universal unital C*-algebra $F$ freely
generated by the family $(z_g)_{g \in G}$ with $\| z_g \| \leq 2$
for each $g \in G.$
(In~\cite{Lr}, the algebra $F$ is not required to be unital, but all
our C*-algebras are unital, so no harm is done by requiring it.)
A set $R$ of relations is simply a subset of $F$ such that, if
$I \subset F$ is the ideal generated by $R,$ then
$\| z_g + I \| \leq 1$ for all $g \in G.$
The universal unital C*-algebra $C^* (G, R)$ is $F / I.$

By Lemma~2.2.5 of~\cite{ELP} (stated in a slightly different setup, but
the difference is unimportant), we may always take $R$ to be finite.
When $G$ and $R$ are finite, we say that $(G, R)$ is a finite
presentation, and that $C^* (G, R)$ is finitely presented.
See Definition~13.2.2 of~\cite{Lr}.

We construct,
by induction on $n,$ unital C*-algebras $B_n \in {\mathcal{C}}$
given by
finite presentations $B_n \cong C^* (G_n, R_n),$
unital homomorphisms
$\ps_{m, n} \colon B_m \to B_n$ for $m \leq n,$ positive integers
$r (n),$ and unital homomorphisms $\bt_n \colon B_n \to A,$ with the
following properties:
\begin{enumerate}
\item\label{Ind1}
$\ps_{m, n} \circ \ps_{l, m} = \ps_{l, n}$ for $l \leq m \leq n.$
\item\label{Ind2}
$\big\| \bt_n \circ \ps_{m, n} (b)
 - \bt_{n - 1} \circ \ps_{m, \, n - 1} (b) \big\| < 2^{-n}$
for $m \leq n - 1$ and $b \in S_{r (n - 1)} (G_m).$
\item\label{Ind3}
$\| \ps_{n - 1, \, n} (b) \| < \| \bt_{n - 1} (b) \| + 2^{- n + 1}$ for
$b \in S_{r (n - 1)} (G_{n - 1}).$
\item\label{Ind4}
$\dist \big(x_k, \, \bt_n (S_{r (n)} (G_n)) \big) < 2^{-n}$
for $0 \leq k \leq n.$
\item\label{Ind5}
$\dist \big( \ps_{m, n} (b), \, S_{r (n)} (G_n) \big) < 2^{-n}$
for $m < n$ and $b \in S_{r (n - 1)} (G_m).$
\item\label{Ind6}
$r (0) < r (1) < r (2) < \cdots$
\end{enumerate}

Suppose the objects above have been constructed.
Let $B = \dirlim B_n,$ with homomorphisms $\ps_n \colon B_n \to B.$
If $b \in S_{r (n)} (G_m)$ and $n \geq m,$ then $k \geq n$ implies
\[
\big\| \bt_{k + 1} \circ \ps_{m, \, k + 1} (b)
                - \bt_k \circ \ps_{m, k} (b) \big\|
< 2^{- k - 1}
\]
by Condition~(\ref{Ind2}), so
$\lim_{k \to \I} \bt_k \circ \ps_{m, k} (b)$ exists.
Since $\bigcup_{n \geq m} S_{r (n)} (G_m)$ is dense in $B_m,$ it
follows that
$\ph_m (b) = \lim_{k \to \I} \bt_k \circ \ps_{m, k} (b)$ exists for
all $b \in B_m,$ and defines a unital homomorphism
$\ph_m \colon B_m \to A$
such that $b \in S_{r (m)} (G_m)$ implies
$\| \ph_m (b) - \bt_m (b) \| < 2^{-m}.$
Moreover, $\ph_n \circ \ps_{m, n} = \ph_m$ whenever $m \leq n,$ so
we get a homomorphism $\ph \colon B \to A$ such that
$\ph \circ \ps_m = \ph_m$ for all $m.$

We prove that $\ph$ is injective.
It suffices to prove that if $n \geq m$ and $b \in S_{r (n - 1)} (G_m),$
then
$\| \ph \circ \ps_m (b) \| > \| \ps_m (b) \| - 2^{-n + 2}.$
By Condition~(\ref{Ind5}), there is $c \in S_{r (n)} (G_n) \subset B_n$
such that $\| \ps_{m, n} (b) - c \| < 2^{-n}.$
We have seen that $\| \ph_n (c) - \bt_n (c) \| < 2^{-n}.$
Therefore, using Condition~(\ref{Ind3}) at the fourth step, we get
\begin{align*}
\| \ps_m (b) \|
& - \| \ph \circ \ps_m (b) \|
   \\
& = \| \ps_n \circ \ps_{m, n} (b) \| - \| \ph_n \circ \ps_{m, n} (b) \|
   \\
& \leq \| \ps_n (c) \| - \| \ph_n (c) \| + 2 \| c - \ps_{m, n} (b) \|
   \\
& \leq \| \ps_{n, \, n + 1} (c) \|
           - \| \bt_n (c) \| + \| \ph_n (c) - \bt_n (c) \|
           + 2 \| c - \ps_{m, n} (b) \|
   \\
& < 2^{-n} + 2^{-n} + 2 \cdot 2^{-n} = 2^{-n + 2}.
\end{align*}
This completes the proof of injectivity.

For surjectivity,
we prove that $\dist (x_m, \, \ph (B)) < 2^{-n + 1}$
for any $n \geq m.$
Use Condition~(\ref{Ind4}) to choose $b \in S_{r (n)} (G_n)$ such that
$\| \bt_n (b) - x_m \| < 2^{-n}.$
Then $\ps_n (b) \in B$ and
\begin{align*}
\| \ph (\ps_n (b)) - x_m \|
 & = \| \ph_n (b) - x_m \|
   \\
 & \leq \| \ph_n (b) - \bt_n (b) \| + \| \bt_n (b) - x_m \|
  < 2^{-n} + 2^{-n} = 2^{-n + 1},
\end{align*}
as desired.

It follows that $A \cong \dirlim B_n,$ and is thus a unital
approximate ${\mathcal{C}}$-algebra.

We now carry out the construction.

To start the induction, we use the hypothesis to find a C*-algebra
$B_0 \in {\mathcal{C}}$
and a unital homomorphism $\bt_0 \colon B_0 \to A$
such that $\dist (x_0, \, \bt_0 (B_0)) < 1.$
Since ${\mathcal{C}}$ is flexible,
and using Lemma~2.2.5 of~\cite{ELP} as described above,
the algebra $B_0$ has a
finite presentation $B_0 \cong C^* (G_0, R_0).$
Since $\bigcup_{r = 0}^{\I} S_r (G_0)$ is dense in $B_0,$ there exists
$r (0)$ such that
$\dist \big( x_0, \, \bt_0 \big( S_{r(0)} (G_0) \big) \big) < 1.$

Now suppose all objects have been constructed through stage $n$ of
the induction.

Since $G_n$ generates $B_n,$ there is $\dt_0 > 0$ such that whenever
$D$ is a unital C*-algebra and $\mu_1, \mu_2 \colon B_n \to D$ are
unital homomorphisms such that $\| \mu_1 (g) - \mu_2 (g) \| < \dt_0$
for all $g \in G_n,$ then $\| \mu_1 (b) - \mu_2 (b) \| < 2^{-n}$
for all
$b \in \bigcup_{m \leq n} \ps_{m, n} \big( S_{r (n)} (G_m) \big).$

Using the notation in the discussion at the beginning of the proof,
recall (Definition~13.2.2 of~\cite{Lr})
that a $\dt$-representation of $(G_n, R_n)$ in a C*-algebra $D$
is family $(a_g)_{g \in G_n}$ in $D$ such that $\| a_g \| \leq 2$
for all $g \in G_n,$
and the homomorphism $\kp \colon F \to D$ such that
$\kp (z_g) = a_g$ for all $g \in G_n$ satisfies
$\| \kp (b) \| \leq \dt$
for all $b \in R.$
Since $B_n \in {\mathcal{C}},$ it is semiprojective,
and Theorem~14.1.4 of~\cite{Lr}
implies that the set $R_n$ of relations is stable.
Therefore there exists $\dt_1 > 0$ such that for every
$\dt_1$-representation $(a_g)_{g \in G_n}$ of $(G_n, R_n)$ in a
C*-algebra $D,$
there is a
representation $(b_g)_{g \in G_n}$ of $(G_n, R_n)$ in $D$ such that
$\| b_g - a_g \| < \frac{1}{2} \dt_0$ for all $g \in G_n.$
Then use finiteness of $R_n$ to choose $\dt_2 > 0$ such that
$\dt_2 \leq
   \min \left(2^{- n - 1}, \, \frac{1}{2} \dt_1, \, \dt_0 \right)$
and so
small that whenever $(c_g)_{g \in G_n}$ is a representation of the
relations $R_n$ in a C*-algebra $D,$ and $\| c_g - d_g \| < \dt_2$
for all $g \in G,$ then $(d_g)_{g \in G}$ is a
$\frac{1}{2} \dt_1$-representation of $(G_n, R_n).$

Define $T = \{ x_0, \ldots, x_{n + 1} \} \cup \bt_n (G_n),$ which is a
finite subset of $A.$
By hypothesis, there exists a C*-algebra $C \in {\mathcal{C}}$
and a unital
homomorphism $\gm \colon C \to A$ such that
$\dist (x, \gm (C)) < \frac{1}{2} \dt_2$ for all $x \in T.$
For each $g \in G_n$ choose $c_g \in C$ such that
$\| \gm (c_g) - \bt_n (g) \| < \frac{1}{2} \dt_2.$
By choice of $\dt_2,$ the elements $\gm (c_g)$ define a
$\frac{1}{2} \dt_1$-representation of $(G_n, R_n).$
(Note that
$\| \gm (c_g) \| < \| \bt_n (g) \| + \frac{1}{2} \dt_2
           \leq 1 + \frac{1}{2} \dt_2 < 2$
for $g \in G_n.$)
Let $\kp \colon F \to A$ be the associated homomorphism, so
$\| \kp (d) \| \leq \frac{1}{2} \dt_1$ for all $d \in R_n.$

By Condition~(\ref{D:Flex:2}) of Definition~\ref{D:Flex}, there
exists an increasing sequence $I_0 \subset I_1 \subset \cdots$ of
ideals in $C$ such that
${\overline{\bigcup_{l = 0}^{\I} I_l}} = {\mathrm{Ker}} (\gm)$ and
$C / I_l \in {\mathcal{C}}$ for every $l.$
Let $\gm = \gm_l \circ \pi_l$ be the factorization of $\gm$ through
$C / I_l.$
Then there exists $l_0$ such that $\| \pi_l (c_g) \| < 2$ for all
$l \geq l_0$ and $g \in G_n.$
Thus, for all $l \geq l_0$ there is a homomorphism
$\kp_l \colon F \to C / I_l$ such that $\kp_l (z_g) = \pi_l (c_g)$ for
$g \in G_n,$ and $\gm_l \circ \kp_l = \kp.$
Since $R_n$ is finite, we may choose $l_1 \geq l_0$ so large that
$l \geq l_1$ implies
$\| \kp_l (d) \| \leq \| \kp (d) \| + \frac{1}{2} \dt_1 \leq \dt_1$
for all $d \in R_n.$
Then $\kp_l$ is a $\dt_1$-representation of $(G_n, R_n).$
By the choice of $\dt_1,$ there is therefore a homomorphism
$\ld \colon B_n \to C / I_{l_1}$ such that
$\| \ld (g) - \kp_l (z_g) \| < \frac{1}{2} \dt_0$ for all $g \in G_n.$
Then for $g \in G_n,$ using $\kp_{l_1} (z_g) = \pi_{l_1} (c_g),$ we
have
\[
\| \gm_{l_1} (\ld (g))  - \bt_n (g) \|
\leq \| \ld (g) - \kp_{l_1} (z_g) \| + \| \gm (c_g) - \bt_n (g) \|
< \tfrac{1}{2} \dt_0 + \tfrac{1}{2} \dt_2
\leq \dt_0.
\]
Now
$\| \gm_{l_1} (\ld (b)) - \bt_n (b) \| < 2^{-n}$ for all
$b \in \bigcup_{m \leq n} \ps_{m,n} \big(S_{r (n)} (G_m) \big)$
by the choice of~$\dt_0.$

For $l \geq l_1,$ let $\mu_l$ be the quotient map
$\mu_l \colon C / I_{l_1} \to C / I_l.$
Then for every $b \in S_{r(n)} (G_n),$ we have
\[
\lim_{l \to \I} \| \mu_l ( \ld (b)) \|
 = \| \gm_{l_1} ( \ld (b)) \|
 < \| \bt_n (b) \| + 2^{-n}.
\]
Therefore there exists $l_2 \geq l_1$ such that
$\| \mu_{l_2} (\ld (b)) \| < \| \bt_n (b) \| + 2^{-n}$ for all
$b \in S_{r (n)} (G_n).$

We now define $B_{n + 1} = C / I_{l_2},$
$\bt_{n +1} = \gm_{l_2},$ and
$\ps_{n, \, n + 1} = \mu_{l_2} \circ \ld.$
For $m < n$ set
$\ps_{m, \, n + 1} = \ps_{n, \, n + 1} \circ \ps_{m, n}.$
Condition~(\ref{Ind1}) is immediate.
Condition~(\ref{Ind2}) follows from the equation
$\gm_{l_1} \circ \ld = \gm_{l_2} \circ \mu_{l_2} \circ \ld$ and the
estimate, proved in slightly different notation above,
$\| \gm_{l_1} \circ \ld \circ \ps_{m, n} (b)
        - \bt_n \circ \ps_{m, n} (b) \| < 2^{-n}$
for $m \leq n$ and $b \in S_{r (n)} (G_m).$
Condition~(\ref{Ind3}) has just been proved.
Since $B_{n + 1} \in {\mathcal{C}},$ there is a finite presentation
$B_{n + 1} \cong C^* (G_{n + 1}, \, R_{n + 1})$ in the sense of
Definition~13.2.2 of~\cite{Lr}.
Then $\bigcup_{r = 0}^{\I} S_r (G_{n + 1})$ is dense in $B_{n + 1}.$
For $k \leq n + 1,$ we have
\[
\dist \big(x_k, \, \bt_{n + 1} (B_{n + 1}) \big)
 = \dist (x_k, \, \gm (C))
 < \dt_2
 \leq 2^{- n - 1},
\]
so it is easy to choose $r (n + 1) > r (n)$ and also so large that
Conditions~(\ref{Ind4}) and~(\ref{Ind5}) are satisfied.
This completes the induction, and thus the proof of the proposition.
\end{proof}

\section{Examples of flexible classes}\label{Sec:ExFlex}

\indent
In this section, we prove that many of the popular classes
of semiprojective building blocks are flexible.

\begin{exa}\label{E:AF}
Let ${\mathcal{C}} = \{ \C \}.$
Then ${\mathcal{C}}$ is flexible.
The finite saturation of ${\mathcal{C}}$
consists of all finite dimensional C*-algebras.
The unital approximate ${\mathcal{C}}$-algebras
are exactly the unital
AF~algebras.
\end{exa}

\begin{exa}\label{E:AI}
Let ${\mathcal{C}} = \{ \C, \, C ([0, 1]) \}.$
We claim that ${\mathcal{C}}$ is flexible.

We verify Condition~(\ref{D:Flex:1}) of Definition~\ref{D:Flex}.
Both $\C$ and $C ([0, 1])$ are well known to be finitely generated
and semiprojective.
Thus, for any $k,$ the algebras $M_k (\C)$ and $M_k ( C ([0, 1]))$
are finitely generated (easy) and semiprojective
(Theorem 14.2.1 of~\cite{Lr}).
Now suppose $A \in {\mathcal{C}}$ and $p \in M_n (A)$ is a projection.
Then there is $k$ such that $p M_n (A) p \cong M_k (A),$
which we just dealt with.

Now we verify Condition~(\ref{D:Flex:2}) of Definition~\ref{D:Flex}.
For $A = \C$ this condition is trivial.
So let $I \subset C ([0, 1])$ be an ideal.
Write
\[
I = \{ f \in C ([0, 1]) \colon f |_Y = 0 \}
\]
for a closed set $Y \subset [0, 1].$
Since $[0, 1] \setminus Y$ is the disjoint union of at most countably
many open intervals, we can write $Y = \bigcap_{n = 0}^{\I} Y_n$ for
closed sets $Y_0 \supset Y_1 \supset \cdots$ such that each
$[0, 1] \setminus Y_n$ is a finite union of open intervals.
Take
\[
I_n = \{ f \in C ([0, 1]) \colon f |_{Y_n} = 0 \}.
\]
Then $C ([0, 1]) / I_n \cong C (Y_n)$ is a finite direct sum of
C*-algebras isomorphic to $\C$ or $C ([0, 1]).$
This completes the proof that ${\mathcal{C}}$ is flexible.

The unital approximate ${\mathcal{C}}$-algebras
are the unital AI~algebras.
Since we have seen several definitions of an AI~algebra,
we formally take this as the definition in this paper.
\end{exa}

\begin{exa}\label{E:AT}
Let ${\mathcal{C}} = \{ \C, \, C ([0,1]), \, C (S^1) \}.$
Then ${\mathcal{C}}$ is flexible
by the same argument as in Example~\ref{E:AI}.
The unital approximate ${\mathcal{C}}$-algebras
are the unital AT~algebras.
Again, we take this as the definition of a unital AT~algebra.
\end{exa}

\begin{exa}\label{E:AD}
For $n \in \N,$ let $D_n$ be the dimension drop interval,
\[
D_n = \{ f \in C ([0, 1], \, M_n ) \colon
         {\mbox{$f(0)$ and $f(1)$ are scalars}} \}.
\]
Further set
\[
C_n = \{ f \in C ([0, 1], \, M_n )
        \colon {\mbox{$f(0)$ is a scalar}} \}.
\]
Let
\[
{\mathcal{C}}
 = \{ \C, \, C ([0, 1]) \} \cup \{ C_n, D_n \colon n \geq 2 \}.
\]

This collection satisfies Condition~(\ref{D:Flex:1}) of
Definition~\ref{D:Flex} because it is contained in the class of
Example~\ref{E:A1D} below.
The argument in Example~\ref{E:AI} for Condition~(\ref{D:Flex:2})
of Definition~\ref{D:Flex} applies directly,
so ${\mathcal{C}}$ is flexible.

We define a unital AD~algebra (``approximate dimension drop algebra'')
to be a unital approximate ${\mathcal{C}}$-algebra.
\end{exa}

One should note that
the collection ${\mathcal{C}}$ in Example~\ref{E:AD}
would not be flexible without the algebras $C_n,$ since
Condition~(\ref{D:Flex:2}) of Definition~\ref{D:Flex} would fail.

We now consider one dimensional noncommutative CW~complexes.
This class requires some preliminary work.
We recall Definition~2.4.1 of~\cite{ELP}
(which is a special case of Definition~11.2 of~\cite{Pd3}):

\begin{dfn}\label{D:NCCW}
A one dimensional noncommutative CW~complex is a \ca\  $A$ of the form
\[
A = \big\{ (a, f) \in F_0 \oplus C ([0, 1], \, F_1)
 \colon {\mbox{$f (0) = \ph_0 (a)$ and $f (1) = \ph_1 (a)$}} \big\},
\]
where $F_0$ and $F_1$ are finite dimensional C*-algebras
and $\ph_0, \ph_1 \colon F_0 \to F_1$ are unital homomorphisms.
\end{dfn}

\begin{lem}\label{L:QNCCW}
Let the notation be as in Definition~\ref{D:NCCW}.
Write $F_1 = \bigoplus_{k = 1}^n M_{d (k)},$
and write $(a, f) \in A$ as
$(a, f_1, f_2, \ldots, f_n)$ with $f_k \in C ([0, 1], \, M_{d (k)}).$
For each $k,$ let $\af_{j, k}$ and $\bt_{j, k},$
for $1 \leq j \leq r (k),$ satisfy
\[
0 \leq \af_{1, k} < \bt_{1, k} < \af_{2, k} < \bt_{2, k}
   < \cdots  < \af_{r (k), \, k} < \bt_{r (k), \, k} \leq 1.
\]
Let $J \S A$ be the ideal consisting of all $(a, f) \in A$
such that $f_k |_{[\af_{j, k}, \, \bt_{j, k}]} = 0$
for $1 \leq k \leq n$ and $1 \leq j \leq r (k).$
Then $A / J$ is a one dimensional noncommutative CW~complex
of the form
\[
A / J = \{ (b, g) \in G_0 \oplus C ([0, 1], \, G_1)
 \colon {\mbox{$g (0) = \ps_0 (b)$ and $g (1) = \ps_1 (b)$}} \},
\]
with $G_0$ and $G_1$ finite dimensional
and $\ps_0, \ps_1 \colon G_0 \to G_1$ unital.
Moreover, there is an expression of this type such that
if $(a, f) \in A$ and $(a, f) + J = (b, g),$
then
\[
\| g \| \leq
  \max_{1 \leq k \leq n} \max_{1 \leq j \leq r (k)}
   \big\| f_k |_{[\af_{j, k}, \, \bt_{j, k}]} \big\|.
\]
\end{lem}

\begin{proof}
Let $\pi \colon A \to A / J$ be the quotient map.

Define $G_1^{(k)} = M_{d (k)}^{r (k)},$ the direct sum of
$r (k)$ copies of $M_{d (k)}.$
Define $G_1 = \bigoplus_{k = 1}^n G_1^{(k)}.$
If $\af_{1, k} = 0$ set
$G_0^{(0, k)} = M_{d (k)}^{r (k) - 1}$;
otherwise set $G_0^{(0, k)} = G_1^{(k)}.$
Similarly, set
$G_0^{(1, k)} = M_{d (k)}^{r (k) - 1}$
if $\bt_{r (k), \, k} = 1$ and
set $G_0^{(1, k)} = G_1^{(k)}$ otherwise.
Further
set
\[
G_0 = F_0 \oplus
 \bigoplus_{k = 1}^n \big( G_0^{(0, k)} \oplus G_0^{(1, k)} \big).
\]
Define $\ps_{0, k} \colon G_0 \to G_1^{(k)}$ as follows.
Let
\[
x = (a, \, b_{0, 1}, \, b_{1, 1}, \, b_{0, 2}, \, b_{1, 2}, \,
  \ldots, \, b_{0, n}, \, b_{1, n})
 \in G_0.
\]
Write $\ph_0 (a)_k$ for the $k$th coordinate of
$\ph_0 (a) \in \bigoplus_{k = 1}^n M_{d (k)}.$
Then set
\[
\ps_{0, k} (x) = \left\{ \begin{array}{ll}
     (\ph_0 (a)_k, \, b_{0, k})   & \hspace{3em}  \af_{1, k} = 0 \\
     b_{0, k}                     & \hspace{3em}  \af_{1, k} > 0
    \end{array} \right.
\]
and
\[
\ps_{1, k} (x) = \left\{ \begin{array}{ll}
     (b_{0, k}, \, \ph_1 (a)_k)   & \hspace{3em}  \bt_{1, k} = 1 \\
     b_{0, k}                     & \hspace{3em}  \bt_{1, k} < 1.
    \end{array} \right.
\]
Combine these in the obvious way to make unital \hm s
$\ps_0, \ps_1 \colon G_0 \to G_1.$
Set
\[
B = \{ (b, g) \in G_0 \oplus C ([0, 1], \, G_1)
 \colon {\mbox{$g (0) = \ps_0 (b)$ and $g (1) = \ps_1 (b)$}} \}.
\]

We show that $A / J \cong B.$
Define $\rh_0 \colon A \to G_0$ by
\[
\rh_0 (a, f)
 = (a, \, b_{0, 1}, \, b_{1, 1}, \, b_{0, 2}, \, b_{1, 2}, \,
  \ldots, \, b_{0, n}, \, b_{1, n}),
\]
with
\[
b_{0, k} = \left\{ \begin{array}{ll}
     \big( f_k (\af_{j, k}) \big)_{2 \leq j \leq r (k)}
                    & \hspace{3em}  \af_{1, k} = 0  \\
     \big( f_k (\af_{j, k}) \big)_{1 \leq j \leq r (k)}
                    & \hspace{3em}  \af_{1, k} > 0
    \end{array} \right.
\]
and
\[
b_{1, k} = \left\{ \begin{array}{ll}
     \big( f_k (\bt_{j, k}) \big)_{1 \leq j \leq r (k) - 1}
                         & \hspace{3em}  \bt_{1, k} = 1  \\
     \big( f_k (\bt_{j, k}) \big)_{1 \leq j \leq r (k)}
                         & \hspace{3em}  \bt_{1, k} < 1.
    \end{array} \right.
\]
For each $j$ and $k,$ choose an orientation preserving homeomorphism
$h_{j, k} \colon [0, 1] \to [\af_{j, k}, \, \bt_{j, k}].$
Then define $\rh_1 \colon A \to C ([0, 1], G_1)$ by
\[
\rh_1 (a, f)
 = \left( \big( f_1 \circ h_{j, 1} \big)_{1 \leq j \leq r (1)},
            \, \big( f_2 \circ h_{j, 2} \big)_{1 \leq j \leq r (2)},
            \, \ldots,
            \, \big( f_n \circ h_{j, n} \big)_{1 \leq j \leq r (n)}
            \right),
\]
and define $\rh \colon A \to B$ by
$\rh (a, f) = (\rh_0 (a, f), \, \rh_1 (a, f) ).$
One checks easily that $\rh (a, f)$ really is in $B$
and that $\rh (a, f) = 0$ \ifo\  $f \in J.$
It remains to show that $\rh$ is surjective.
Let $(b, g) \in B.$
Let $a$ be the projection of $b$ to the summand $F_0$ of $G_0.$
Let $1 \leq k \leq n.$
Set
\[
L_k
 = \{ 0, 1 \} \cup \bigcup_{j = 1}^{r (k)} [\af_{j, k}, \, \bt_{j, k}].
\]
We claim that there exists a unique function
$f_k^{(0)} \in C (L_k, M_{d (k)})$ such that
$f_k^{(0)} (0) = \ph_0 (a)_k,$ $f_k^{(0)} (1) = \ph_1 (a)_k,$
and the functions $f_k^{(0)} \circ h_{j, k}$
are equal to the corresponding components of $g.$
This is immediate unless $\af_{1, k} = 0$ or $\bt_{r (k), k} = 1,$
in which case $f_k^{(0)} (0)$ or $f_k^{(0)} (1)$ must satisfy
two conditions.
In these cases, however, the requirements
$g (0) = \ps_0 (b)$ and $g (1) = \ps_1 (b)$ ensure that the
two conditions agree.
Now extend $f_k^{(0)}$ to a continuous function
$f_k \in C ([0, 1], M_{d (k)}),$
and observe that $(a, f_1, f_2, \ldots, f_n) \in A$
and satisfies $\rh (a, f_1, f_2, \ldots, f_n) = (b, g).$
This completes the proof,
except for the norm estimate, which is clear from the construction.
\end{proof}

\begin{lem}\label{L:NCCWAId}
Let $A$ be a one dimensional noncommutative CW~complex,
and let $I$ be an ideal in $A.$
Then there exists an
increasing sequence $I_0 \subset I_1 \subset \cdots$ of ideals in
$A$ such that ${\overline{\bigcup_{n = 0}^{\I} I_n}} = I,$ and such
that for every $n$ the C*-algebra $A / I_n$
is a one dimensional noncommutative CW~complex.
\end{lem}

\begin{proof}
Let the notation be as in Definition~\ref{D:NCCW}.
We claim that it suffices to consider the case in which
$(\ph_0, \ph_1) \colon F_0 \to F_1 \oplus F_1$ is injective.
If this map is not injective, let $L$ be its kernel,
and write $A = L \oplus B$ where $B$
is the one dimensional noncommutative CW~complex defined
using $F_0 / L$ and the induced maps $F_0 / L \to F_1$
in place of $F_0,$ $\ph_0,$ and $\ph_1.$
It suffices to consider the summands separately,
so the claim follows.

We first claim that, under the hypotheses of the lemma,
and assuming the injectivity condition of the previous paragraph,
if $x \in I$ then
for every $\ep > 0$ there is an ideal $J \S I$
such that $A / J$ is a one dimensional noncommutative CW~complex
satisfying the injectivity condition of the previous paragraph
and such that $\| x + J \| < \ep.$
Write $F_1 = \bigoplus_{k = 1}^n M_{d (k)},$
and write $x = (a, f_1, f_2, \ldots, f_n)$
with $a \in F_0$ and $f_k \in C ([0, 1], \, M_{d (k)}).$
Let
\[
U_k = \big\{ t \in [0, 1] \colon \| f_k (t) \| > 0 \big\}
\andeqn
L_k = \big\{ t \in [0, 1]
   \colon \| f_k (t) \| \geq \tfrac{1}{2} \ep \big\}.
\]
Since $U_k$ is a countable disjoint union of open intervals in $[0, 1]$
(we include intervals of the form $[0, s)$ and $(s, 1]$),
we may write $U_k = \bigcup_{m = 1}^{\infty} V_{m, k},$
where each $V_{m, k}$ is a finite union of open intervals,
the closures of which are contained
in distinct intervals in the expression of $U_k$
as a disjoint union of open intervals.
By compactness, one of these sets, call it $W_k,$ contains $L_k.$
Let $J_0$ be the ideal in $A$ consisting of all
$(b, g_1, g_2, \ldots, g_n) \in A$
such that $g_k |_{[0, 1] \setminus W_k} = 0$ for $1 \leq k \leq n.$
The sets $[0, 1] \setminus W_k$ have the form in the hypotheses
of Lemma~\ref{L:QNCCW},
so $A / J_0$ is a one dimensional noncommutative CW~complex
of the form
\[
A / J_0 = \{ (b_0, g) \in G_0 \oplus C ([0, 1], \, G_1)
 \colon {\mbox{$g (0) = \ps_0 (b)$ and $g (1) = \ps_1 (b)$}} \},
\]
such that if we write $x + J_0 = (b, g),$
then $\| g \| < \ep.$
Set $M = \Ker (\ps_0) \cap \Ker (\ps_1).$
Define a new one dimensional noncommutative CW~complex
$B$ by replacing $G_0$ with $G_0 / M,$
and making the obvious replacements for $\ps_0$ and $\ps_1.$
The map $G_0 / M \to G_1 \oplus G_1$ they induce is injective.
Thus, letting $J$ be the kernel of $A \to B,$
we have $\| x + J \| \leq \| g \| < \ep.$
This completes the proof of the claim.

Now we prove the lemma.
Let $(x_n)_{n \in \N}$ be a sequence in $I$ such that
every tail $(x_n)_{n \geq N}$ is dense.
We apply the claim to $I \subset A$
with $x_1$ in place of $x$ and with $\ep = 1,$
obtaining an ideal $I_1 \subset I.$
Let $\pi_1 \colon A \to A / I_1$ be the quotient map.
Now apply the claim to $I / I_1 \subset A / I_1$
with $\pi_1 (x_2)$ in place of $x$ and with $\ep = \frac{1}{2},$
obtaining an ideal $J_2 \subset I / I_1.$
Take $I_2 = \pi_1^{-1} (J_2),$
and let $\pi_2 \colon A \to A / I_2$ be the quotient map.
Then apply the claim to $I / I_2 \subset A / I_2$
with $\pi_2 (x_3)$ in place of $x$ and with $\ep = \frac{1}{3},$
obtaining an ideal $J_3 \subset I / I_2,$
and take $I_3 = \pi_2^{-1} (J_3),$
etc.
For every $N,$ every element of the dense sequence
$(x_n)_{n \geq N}$ in $I$ can be approximated within $1 / N$
by elements of $\bigcup_{n = 1}^{\infty} I_n,$
so that $\bigcup_{n = 1}^{\infty} I_n$ is dense in $I.$
\end{proof}

\begin{exa}\label{E:A1D}
Let ${\mathcal{C}}$
be the collection of all \ca s
isomorphic to one dimensional noncommutative
CW~complexes as in Definition~\ref{D:NCCW}.
We claim that ${\mathcal{C}}$ is flexible.
One easily checks ${\mathcal{C}}$ is finitely saturated.
Algebras in ${\mathcal{C}}$
are finitely generated by Lemma~2.4.3 of~\cite{ELP},
and semiprojective by Theorem~6.2.2 of~\cite{ELP}.
Thus, Condition~(\ref{D:Flex:1}) of Definition~\ref{D:Flex}
is satisfied.
Condition~(\ref{D:Flex:2}) of Definition~\ref{D:Flex}
follows from Lemma~\ref{L:NCCWAId}.
\end{exa}

\begin{exa}\label{E:AH}
For $d \in \N,$ take ${\mathcal{C}}_d$
to be the collection of all
C*-algebras $C (X)$ for $X$ a compact metric space with covering
dimension at most $d.$
Then a separable unital C*-algebra $A$ is an AH~algebra with no
dimension growth if and only if it is a unital approximate
${\mathcal{C}}_d$-algebra for some $d.$

The class ${\mathcal{C}}_d$ is not flexible.
(For $d = 0,$ it contains the algebra of continuous functions on the
Cantor set, which is not semiprojective.)
Even if one restricts to finite complexes, one will not get a
flexible class for $d \geq 2.$
\end{exa}

\begin{exa}\label{R:OtherFlex}
Let $T$ be the Toeplitz algebra, the \ca\  of the unilateral shift.
Then the class
$\{ \C, \, C ([0, 1]), \, C (S^1), \, T \}$
is flexible.
The verification
of Condition~(\ref{D:Flex:2}) of Definition~\ref{D:Flex}
for $T$ reduces to that for $C (S^1).$
For Condition~(\ref{D:Flex:1}),
let $p \in M_n (T)$ be a nonzero \pj.
Let $\pi \colon M_n (T) \to M_n (C (S^1))$ be the quotient map.
If $\pi (p) = 0,$ then $p M_n (T) p$ is finite dimensional,
hence obviously finitely generated and projective.
Otherwise, there is a short exact sequence
\[
0 \longrightarrow p K p \longrightarrow p M_n (T) p
    \longrightarrow \pi (p) M_n (C (S^1))  \pi (p) \longrightarrow 0.
\]
Clearly $\pi (p)$ is full, so $p$ is full.
Now $T$ is semiprojective by Corollary~2.22 of~\cite{Bl2},
so $M_n (T)$ is semiprojective by Theorem~14.2.2 of~\cite{Lr},
whence
$p M_n (T) p$ is semiprojective by Proposition~2.27 of~\cite{Bl2}.
Moreover, $\pi (p) M_n (C (S^1))  \pi (p)$ is finitely generated
as in Example~\ref{E:AT},
and $p K p \cong K$ is finitely generated by Theorem~8 of~\cite{OZ},
so $p M_n (T) p$ is finitely generated.
\end{exa}

We now consider examples using \pisca s.
The following lemma is useful.
The conclusion can be improved to ``singly generated'';
see Example~(4) on page~138 of~\cite{Ng}.
However, no proof is given there,
and the version here is easy to prove and suffices for our purposes.

\begin{lem}\label{L:PIFG}
Every separable unital \pisca\  is finitely generated.
\end{lem}

\begin{proof}
Let $A$ be a separable unital \pisca.
Choose nonzero \pj s $p_1, p_2, p_3 \in A$ such that
$p_1 + p_2 + p_3 = 1.$
There exist \mops\  $e_1, e_2, \ldots \leq p_1,$ each equivalent
to $1 - p_1,$
and the \hsa\  $B$ of $A$ generated by
$\{ 1 - p_1 \} \cup \{ e_1, e_2, \ldots \}$
is isomorphic to $K \otimes (1 - p_1) A (1 - p_1).$
It follows from Theorem~8 of~\cite{OZ} that $B$ is generated by
a single element $a_1.$
Similarly, there exist $a_2, a_3 \in A$ which generate subalgebras
which contain $(1 - p_2) A (1 - p_2)$ and $(1 - p_3) A (1 - p_3).$
In particular, $\{ a_1, a_2, a_3 \}$ generates a subalgebra $A_0$ of $A$
which contains $(1 - p_j) A (1 - p_j)$ for $j \in \{ 1, 2, 3 \}.$
If $j, k \in \{ 1, 2, 3 \}$ are given, and $l$ is equal to
neither $j$ nor $k,$ then $p_j A p_k \subset (1 - p_l) A (1 - p_l).$
So $A_0$ contains $p_j A p_k$ for all $j$ and $k,$
and hence contains $A.$
\end{proof}

\begin{exa}\label{E:Cuntz}
The class of Cuntz algebras
${\mathcal{C}} = \{ {\mathcal{O}}_m \colon 2 \leq m < \infty \}$
is flexible.
For the proof,
${\mathcal{O}}_m$ is semiprojective by Corollary~2.24 of~\cite{Bl2},
so $M_n ({\mathcal{O}}_m)$ is semiprojective
by Theorem~14.2.2 of~\cite{Lr}.
Every corner in such an algebra is full,
hence semiprojective by Proposition~2.27 of~\cite{Bl2}.
It follows from Lemma~\ref{L:PIFG} that all such corners are finitely
generated.

Every (nonempty) subset of ${\mathcal{C}}$ is also flexible.
\end{exa}

\begin{exa}\label{E:OInf}
The class $\{ {\mathcal{O}}_{\infty} \}$ is flexible.
Indeed, this algebra is semiprojective by Theorem~3.2 of~\cite{Bl7},
and the rest of the argument is as in Example~\ref{E:Cuntz}.
\end{exa}

\begin{exa}\label{E:Union}
The union of flexible classes is flexible.
Thus, for example, the class
$\{ \C, \, C ([0, 1]), \, C (S^1), \, {\mathcal{O}}_2 \}$
is flexible,
by Examples \ref{E:AT} and~\ref{E:Cuntz},
and the class $\{ {\mathcal{O}}_m \colon 2 \leq m \leq \infty \}$
is flexible, by Examples \ref{E:Cuntz} and~\ref{E:OInf}.
\end{exa}

\begin{exa}\label{E:Em}
Let $E_m$ be the standard extension of ${\mathcal{O}}_m$ by $K.$
Then the class $\{ \C, {\mathcal{O}}_m, E_m \}$ is flexible.
The algebra ${\mathcal{O}}_m$ is semiprojective by
Corollary~2.24 of~\cite{Bl2}.
The algebra $E_m$ is easily seen to be semiprojective
by using Propositions 2.18 and~2.23 of~\cite{Bl2}.
The proof of Condition~(\ref{D:Flex:1}) of Definition~\ref{D:Flex}
is now the same as for Example~\ref{R:OtherFlex},
but using semiprojectivity and finite generation of
corners and matrix algebras over ${\mathcal{O}}_m$ as in
Example~\ref{E:Cuntz}.
\end{exa}

\section{Structure theorems for crossed products by
    actions with the Rokhlin property}\label{Sec:CrPrdRP}

\indent
In this section,
we prove various results of the form discussed in the introduction,
to the effect that crossed products by finite group actions
with the Rokhlin property preserve various properties of \ca s.
We concentrate in membership in various classes considered
in the Elliott program;
other properties are considered in the next section.
The Rokhlin property is as in Definition~3.1 of~\cite{Iz1};
we reproduce here the equivalent form given in
Definition~1.1 of~\cite{PhtRp1a}.

\begin{dfn}\label{SRPDfn}
Let $A$ be a unital \ca,
and let $\af \colon G \to \Aut (A)$
be an action of a finite group $G$ on $A.$
We say that $\af$ has the
{\emph{Rokhlin property}} if for every finite set
$F \S A$ and every $\ep > 0,$
there are \mops\  $e_g \in A$ for $g \in G$ such that:
\begin{enumerate}
\item\label{SRPDfn:1}
$\| \af_g (e_h) - e_{g h} \| < \ep$ for all $g, h \in G.$
\item\label{SRPDfn:2}
$\| e_g a - a e_g \| < \ep$ for all $g \in G$ and all $a \in F.$
\item\label{SRPDfn:3}
$\sum_{g \in G} e_g = 1.$
\end{enumerate}
\end{dfn}

The proof of the following theorem is contained in
the proof of Theorem~2.2 of~\cite{PhtRp1a}.
We will mostly use the version in Theorem~\ref{T:CrPrdRk} below,
but one of our results requires this version.

\begin{thm}\label{T:CrPrdRk_0}
Let $A$ be a unital \ca,
let $G$ be a finite group,
and let $\af \colon G \to \Aut (A)$ be an
action with the Rokhlin property.
Then for every finite subset $S \S C^* (G, A, \af)$ and every $\ep > 0,$
there are $n,$ a \pj\  $f \in A,$
and a unital \hm\  $\ph \colon M_n \otimes f A f \to C^* (G, A, \af)$
such that $\dist (a, \, \ph (M_n \otimes f A f)) < \ep$
for all $a \in S.$
\end{thm}

\begin{proof}
This is what is actually shown
in the proof of Theorem~2.2 of~\cite{PhtRp1a}.
(We take $f$ to be the \pj\  $e_1$ used there.)
\end{proof}

Theorem~\ref{T:CrPrdRk_0} has the following consequence,
which is usually what is actually wanted.
When $A$ is assumed to be an AH~algebra with no dimension growth,
the conclusion was first obtained by Dawn Archey.

\begin{thm}\label{T:CrPrdRk}
Let ${\mathcal{C}}$ be any class of separable unital C*-algebras
which has approximate quotients (Definition~\ref{D:Flex}).
Let $A$ be a unital local ${\mathcal{C}}$-algebra.
Let $G$ be a finite group,
and let $\af \colon G \to \Aut (A)$ be an
action with the Rokhlin property.
Then $C^* (G, A, \af)$ is again a unital local ${\mathcal{C}}$-algebra.
\end{thm}

\begin{proof}
Parts (\ref{L:ClosULC:2}) and~(\ref{L:ClosULC:5})
of Lemma~\ref{L:ClosULC}
show that the algebras $M_n \otimes f A f$ in Theorem~\ref{T:CrPrdRk_0}
are unital local ${\mathcal{C}}$-algebras.
Therefore Theorem~\ref{T:CrPrdRk_0}
and Lemma~\ref{L:ClosULC}(\ref{L:ClosULC:3a})
imply that $C^* (G, A, \af)$ is a unital local ${\mathcal{C}}$-algebra.
\end{proof}

In Theorem~\ref{T:CrPrdRk_0},
one might hope that one can choose
$f$ so that $f A f$ is stably isomorphic to a
direct summand of $A.$
If so, in Theorem~\ref{T:CrPrdRk}
one could replace having approximate quotients
by a weaker hypothesis.
One would obtain,
for example,
a strengthening of Proposition~\ref{C:RokhUCT} below,
in which nothing need be said about quotients.
Unfortunately, the following example shows that, in general,
it is not possible to choose $f$ to satisfy this condition.
(The example does not, however,
show that the proposed strengthening of Proposition~\ref{C:RokhUCT}
is false.)

\begin{exa}\label{E:NoDS}
Let $B = \bigotimes_{n = 1}^{\I} M_2$ be the $2^{\infty}$~UHF algebra.
We set $B_n = \bigotimes_{k = 1}^{n} M_2,$
and identify $B$ with $\dirlim B_n.$
Let
$\bt_0 = \bigotimes_{n = 1}^{\I} \Ad \left( \left(
  \begin{smallmatrix} 1 & 0 \\ 0 & -1 \end{smallmatrix} \right) \right)
 \in \Aut (B).$
Let $K = K (l^2 (\N)),$
and for $j, k \in \N$ let $f_{j, k}$ denote the standard matrix unit.
Further let $f_n = \sum_{k = 1}^{2 n} f_{k, k}$
and $g_n = \sum_{k = 1}^n f_{2 k, 2 k}.$
In the multiplier algebra $M (K \otimes B)$ define
(with convergence in the strict topology)
\[
g = \lim_{n \to \infty} g_n \otimes 1
  = \sum_{k = 1}^{\infty} f_{2 k, 2 k} \otimes 1
\andeqn
u = \sum_{k = 1}^{\infty}
   ( f_{2 k - 1, \, 2 k} + f_{2 k, \, 2 k - 1} ) \otimes 1.
\]
Then $u$ is unitary, $u^2 = 1,$ and $u g u^* = 1 - g.$
Define $A \subset M (K \otimes B)$ to be the subalgebra
generated by $1,$ $g,$ and $K \otimes B.$
Let $\bt \in \Aut ( M (K \otimes B) )$ be the extension to
$M (K \otimes B)$ of $\id_K \otimes \bt_0,$
and let $\af \in \Aut (A)$ be the restriction to $A$
of $\Ad (u) \circ \bt.$
Note that $(\Ad (u) \circ \bt) (A) = A,$
so that we really do have $\af \in \Aut (A).$
Moreover, since $\Ad (u)$ and $\bt$ are commuting automorphisms
of order~$2,$
we easily see that $\af^2 = \id_A.$

We claim that $\af$ generates an action of $\Z / 2 \Z$
which has the Rokhlin property.
Let $F_n \S A$ be the set
\[
F_n = \{1,  g \} \cup
         \big\{ e_{j, k} \otimes b \colon
            {\mbox{$j, k \leq 2 n$ and $b \in B_n$}} \big\}.
\]
The sets $F_n$ are not finite, but
$F_1 \subset F_2 \subset \cdots$ and
$\bigcup_{n = 1}^{\I} F_n$ generates a dense subalgebra of~$A,$
so it suffices to verify the Rokhlin property using the
sets $F_n$ in place of the finite sets of Definition~\ref{SRPDfn}.

So fix~$n.$
Let $p \in B$ be the \pj\  %
\[
p = 1 \otimes 1 \otimes \cdots \otimes 1
        \otimes
        \frac{1}{2}
          \left( \begin{matrix} 1 & 1 \\ 1 & 1 \end{matrix} \right)
        \otimes 1 \otimes 1 \otimes \cdots,
\]
with the nontrivial entry in the tensor factor in position $n + 1.$
Then $p$ commutes with everything in $B_n$ and $\bt_0 (p) = 1 - p.$
Now set
\[
q = g - g_n \otimes 1 + f_n \otimes p.
\]
Using the relations
\[
u (f_n \otimes 1) = (f_n \otimes 1) u,
\,\,\,\,\,\,
u (g_n \otimes 1) u^* = (f_n - g_n) \otimes 1,
\andeqn
u g u^* = 1 - g,
\]
we get $\af (q) = 1 - q.$
Moreover, one easily checks that
\[
g - g_n \otimes 1 = \sum_{k = n + 1}^{\infty} f_{2 k, 2 k} \otimes 1
\andeqn
f_n \otimes p
\]
commute with every element of $F_n.$
This shows that $\af$ generates an action with the Rokhlin property,
proving the claim.

Now we claim that if, in the definition of the Rokhlin property,
we take $\ep < 1,$
then it is not possible to choose $e_1$ such that
$e_1 A e_1$ is stably isomorphic to a direct summand of~$A.$
Since the \pj\  $f$ in Theorem~\ref{T:CrPrdRk_0}
is just the \pj\  $e_1$ in an application of the Rokhlin property,
with $\ep > 0$ arbitrarily small,
this will complete the proof
that our example has the properties claimed for it.

Since $A$ can't be written as a direct sum in a nontrivial way,
it suffices to rule out
stable isomorphism with $A$ and the zero algebra.
We can disregard the zero algebra,
since $\ep < 1$ implies $e_1 \neq 0.$
Let $\pi \colon A \to \C \oplus \C$ be the quotient map
with kernel $K \otimes B$ and which sends $g$ to $(1, 0).$
Let $\gm \in \Aut (\C \oplus \C)$ exchange the summands.
Then $\gm \circ \pi = \pi \circ \af.$
Since $\| \af (e_1) - (1 - e_1) \| < 1,$
we must have $\| \gm ( \pi (e_1)) - (1 - \pi (e_1)) \| < 1.$
This can only happen if $\pi (e_1)$ is either $(1, 0)$ or $(0, 1).$
In either case,
$\pi (e_1 A e_1) \cong \C,$
so $e_1 A e_1$ has at most one nontrivial ideal
(namely $e_1 (K \otimes B) e_1$).
This is not true of~$A,$
so $e_1 A e_1$ is not stably isomorphic to~$A.$
This completes the proof of the claim.
\end{exa}

\begin{thm}\label{T:CrPrdFlx}
Let ${\mathcal{C}}$
be any flexible class of separable unital C*-algebras
(Definition~\ref{D:Flex}).
Let $A$ be a unital approximate ${\mathcal{C}}$-algebra.
Let $G$ be a finite group,
and let $\af \colon G \to \Aut (A)$ be an
action with the Rokhlin property.
Then $C^* (G, A, \af)$
is again a unital approximate ${\mathcal{C}}$-algebra.
\end{thm}

\begin{proof}
Combine Theorem~\ref{T:CrPrdRk} and Proposition~\ref{P:LCFlex}.
\end{proof}
In particular, we have the following result.

\begin{cor}\label{C:ACR}
Let $A$ be a separable unital C*-algebra,
let $G$ be a finite group,
and let $\af \colon G \to \Aut (A)$ have the Rokhlin property.
\begin{enumerate}
\item\label{C:AIR}
If $A$ is a unital AI~algebra,
as defined in Example~\ref{E:AI},
then $C^* (G, A, \af)$ is a unital
AI~algebra.
\item\label{C:ATR}
If $A$ is a unital AT~algebra,
as defined in Example~\ref{E:AT},
then $C^* (G, A, \af)$ is a unital
AT~algebra.
\item\label{C:ADR}
If $A$ is a unital AD~algebra,
as defined in Example~\ref{E:AD},
then $C^* (G, A, \af)$ is a unital
AD~algebra.
\item\label{C:A1D}
If $A$ is a unital countable direct limit of
one dimensional noncommutative CW~complexes (Definition~\ref{D:NCCW}),
then so is $C^* (G, A, \af).$
\end{enumerate}
\end{cor}

\begin{proof}
The relevant classes are flexible by Examples~\ref{E:AI},
\ref{E:AT}, \ref{E:AD}, and~\ref{E:A1D},
so we may apply Theorem~\ref{T:CrPrdFlx}.
\end{proof}

The same conclusion holds for the other flexible classes
given in Section~\ref{Sec:ExFlex}.

We now want to prove similar results for
simple unital AH~algebras with slow dimension growth
and real rank zero,
and for Kirchberg algebras satisfying the Universal Coefficient Theorem.
These results require preparation. 

\begin{prp}\label{C:RokhUCT}
Let $A$ be a separable nuclear unital C*-algebra,
and suppose that $A / I$
satisfies the Universal Coefficient Theorem
for every ideal $I \S A.$
Let $G$ be a finite group,
and let $\af \colon G \to \Aut (A)$ have the Rokhlin property.
Then $C^* (G, A, \af) / J$ satisfies the Universal Coefficient Theorem
for every ideal $J \S C^* (G, A, \af).$
\end{prp}

\begin{proof}
Take ${\mathcal{C}}$ to be the set of all quotients of $A.$
The finite saturation of ${\mathcal{C}}$ consists of finite
direct sums of unital \ca s which are stably isomorphic to
quotients of~$A.$
All quotients of such algebras are again algebras of the same form.
So $C^* (G, A, \af)$ is a unital local ${\mathcal{C}}$-algebra
by Theorem~\ref{T:CrPrdRk}.
Let $J \S C^* (G, A, \af)$ be an ideal.
Then $C^* (G, A, \af) / J$
is also a unital local ${\mathcal{C}}$-algebra
by Lemma~\ref{L:ClosULC}(\ref{L:ClosULC:3}).
The algebras in the finite saturation of ${\mathcal{C}},$
and their quotients,
are all nuclear and all satisfy the Universal Coefficient Theorem.
It now follows from Theorem~1.1 of~\cite{Dd}
that $C^* (G, A, \af) / J$ satisfies the Universal Coefficient Theorem.
\end{proof}

\begin{rmk}\label{R:KThyRkh}
Similar arguments can be used to show that,
for example,
if $G$ is finite,
$\af \colon G \to \Aut (A)$ has the Rokhlin property,
and $K_0 (A/I) = 0$ for every ideal $I \S A,$
then $K_0 (C^* (G, A, \af) / J) = 0$
for every ideal $J \S C^* (G, A, \af).$
However, at least for simple $A,$
much more is known:
$K_* (A^G) \to K_* (A)$ is injective by Theorem~3.13 of~\cite{Iz1},
and $K_* (A^G) \cong K_* ( C^* (G, A, \af))$
because $A^G$ is isomorphic to a corner in $C^* (G, A, \af)$
and $C^* (G, A, \af)$ is simple.
\end{rmk}

\begin{cor}\label{C:RokhSUCT}
Let $A$ be a simple separable nuclear unital C*-algebra
which satisfies the Universal Coefficient Theorem,
let $G$ be a finite group,
and let $\af \colon G \to \Aut (A)$ have the Rokhlin property.
Then $C^* (G, A, \af)$ satisfies the Universal Coefficient Theorem.
\end{cor}

\begin{thm}\label{C:RokhAH}
Let $A$ be a simple unital AH~algebra with slow dimension growth
and real rank zero,
let $G$ be a finite group,
and let $\af \colon G \to \Aut (A)$ have the Rokhlin property.
Then $C^* (G, A, \af)$ is a simple unital AH~algebra
with slow dimension growth and real rank zero.
\end{thm}

\begin{proof}
The \ca\  $A$ has tracial rank zero by Corollary~6.2.5 of~\cite{Lnbk}.
Therefore $C^* (G, A, \af)$ has tracial rank zero
by Theorem~2.2 of~\cite{PhtRp1a}.
Also, $C^* (G, A, \af)$ is simple
(by Remark~1.4 and Corollary~1.6 of~\cite{PhtRp1a}),
separable, and nuclear,
and satisfies the Universal Coefficient Theorem
(by Corollary~\ref{C:RokhSUCT}).
Therefore the classification theorem of~\cite{Ln15},
in the form given in Proposition~3.7 of~\cite{PhtRp2},
implies that $C^* (G, A, \af)$ is a simple unital AH~algebra
with slow dimension growth and real rank zero.
\end{proof}

Unlike the results in Corollary~\ref{C:ACR},
the proof requires the additional assumptions
of simplicity and real rank zero.

A similar result holds for Kirchberg algebras.

\begin{cor}\label{C:RokhKr}
Let $A$ be a unital Kirchberg algebra
(a simple, separable, nuclear, and purely infinite \ca)
satisfying the Universal Coefficient Theorem,
let $G$ be a finite group,
and let $\af \colon G \to \Aut (A)$ have the Rokhlin property.
Then $C^* (G, A, \af)$ is a unital Kirchberg algebra
satisfying the Universal Coefficient Theorem.
\end{cor}

\begin{proof}
It is well known that the crossed product is nuclear.
The action $\af$ is outer, by Lemma~1.5 of~\cite{PhtRp1a}.
Therefore $C^* (G, A, \af)$ is a unital Kirchberg algebra,
by Theorem~3 of~\cite{Je}.
It satisfies the Universal Coefficient Theorem
by Corollary~\ref{C:RokhSUCT}.
\end{proof}

\section{Further properties preserved by crossed products by
    actions with the Rokhlin property}\label{Sec:CrPrdRP2}

In this section, we consider preservation of real and stable rank
and related properties,
such as approximate divisibility, isometric richness,
and having the order on \pj s be determined by traces.
We also prove that all tracial states on a crossed
product by a finite group action with the Rokhlin property
are induced from invariant tracial states on the original algebra.

\begin{prp}\label{P:RkOfCrPrd}
Let $A$ be a separable unital C*-algebra,
let $G$ be a finite group,
and let $\af \colon G \to \Aut (A)$ have the Rokhlin property.
\begin{enumerate}
\item\label{PR:TSR}
If $A$ has stable rank one, then so does $C^* (G, A, \af).$
\item\label{PR:RR}
If $A$ has real rank zero, then so does $C^* (G, A, \af).$
\end{enumerate}
\end{prp}

\begin{proof}
For~(\ref{PR:TSR}),
let ${\mathcal{C}}$ be the class of all separable unital C*-algebras
with stable rank one.
We claim that ${\mathcal{C}}$ is finitely saturated.
Conditions (\ref{D:FSat:1}) and~(\ref{D:FSat:2}) are obvious,
and Condition~(\ref{D:FSat:3}) is Theorem~3.3 of~\cite{Rf}.
For Condition~(\ref{D:FSat:4}) (corners),
we observe that $p A p$ is stably isomorphic to the ideal
it generates, and use Theorems 3.6 and~4.4 of~\cite{Rf}.
This proves the claim.
It further follows from Theorem~4.3 of~\cite{Rf}
that ${\mathcal{C}}$ is closed under passage to quotients.
Theorem~\ref{T:CrPrdRk} now implies that
for every finite set $F \subset C^* (G, A, \af)$ and every $\ep > 0,$
there is $B \in {\mathcal{C}}$
and a
unital homomorphism $\ph \colon B \to C^* (G, A, \af)$ such that
$\dist (a, \, \ph (B)) < \ep$ for all $a \in F.$
This clearly implies that $C^* (G, A, \af)$ has stable rank one.

For~(\ref{PR:RR}),
the same argument applies,
using the class of separable unital C*-algebras with real rank zero,
and using Theorems 2.5, 2.10, and~3.14 of~\cite{BP}.
\end{proof}

One expects that a weaker condition on the action than the
Rokhlin property should suffice for Proposition~\ref{P:RkOfCrPrd}.
Compare with the results of~\cite{OP1} for crossed products by~$\Z,$
and analogous results of~\cite{Ar} for the case of finite groups
acting on simple \ca s.
However, the proposition fails
if no condition at all is put on the action.
See Example 8.2.1 of~\cite{Bl3} for stable rank one,
and Example~9 of~\cite{El} for real rank zero.

The same proof applies
to finite stable rank in place of stable rank one,
with appropriate changes in the justification of finite saturation.
However, in this case,
no condition on $\af$ is needed.
See Theorem~2.2 and Example~2.1 of~\cite{JOPT}.
For finite real rank, we do not know enough about the real rank
of corners and matrix algebras.

Recall (see~4.1 of~\cite{BP4}) that a unital \ca\  $A$ is said
to be isometrically rich if the set of one-sided invertible
elements is dense in~$A.$
This is a weakening of stable rank one which holds for purely
infinite simple \ca s.
It is stronger than extremal richness,
as defined at the beginning of Section~3 of~\cite{BP2}.

\begin{prp}\label{P:IsoRich}
Let $A$ be an isometrically rich unital  \ca,
let $G$ be a finite group,
and let $\af \colon G \to \Aut (A)$ have the Rokhlin property.
Then $C^* (G, A, \af)$ is isometrically rich.
\end{prp}

For the proof, we can't use results about finitely saturated classes.
The problem is that
the direct sum of two unital isometrically rich \ca s
is generally not isometrically rich.
If $a \in A$ is right invertible but not invertible,
and $b \in B$ is left invertible but not invertible,
then $(a, b) \in A \oplus B$
can't be approximated by one-sided invertible elements.

\begin{proof}[Proof of Proposition~\ref{P:IsoRich}]
Let $a \in C^* (G, A, \af),$
and let $\ep > 0.$
Apply Theorem~\ref{T:CrPrdRk_0} with $\frac{1}{2} \ep$ in place of $\ep$
and with $S = \{ a \},$
obtaining $\rh \colon M_n \otimes f A f \to C^* (G, A, \af)$
and $b \in M_n \otimes f A f$
such that $\| \rh (b) - a \| < \frac{1}{2} \ep.$

We claim that $M_n \otimes f A f$ is isometrically rich.
Proposition~4.5 of~\cite{BP4} provides an
extremally rich primitive \ca~$B$
and a surjective \hm\  $\pi \colon B \to A.$
Choose $b \in B$ with $b \geq 0$ such that $\pi (b) = f.$
Then ${\overline{b B b}}$ is extremally rich
by Theorem~3.5 of~\cite{BP2},
and is still primitive;
moreover, $\pi$ restricts to a surjective \hm\  %
from ${\overline{b B b}}$ to $f A f.$
Tensor with $\id_{M_n}$ to get a surjective \hm\  %
from $M_n \otimes {\overline{b B b}}$ to $M_n \otimes f A f.$
The algebra $M_n \otimes {\overline{b B b}}$ is extremally rich
by Theorem~4.5 of~\cite{BP2},
and is still primitive.
So Proposition~4.5 of~\cite{BP4} implies that
$M_n \otimes f A f$ is isometrically rich,
proving the claim.

Accordingly,
there is a one-sided invertible element $c \in M_n \otimes f A f$
such that $\| c - b \| < \frac{1}{2} \ep.$
Then $\rh (c)$ is one-sided invertible in $\rh (M_n \otimes f A f),$
and hence also in $C^* (G, A, \af).$
Also, $\| \rh (c) - a \| < \ep.$
\end{proof}

\begin{cor}\label{C:ExtR}
Let $A$ be an extremally rich unital prime  \ca,
let $G$ be a finite group,
and let $\af \colon G \to \Aut (A)$ have the Rokhlin property.
Then $C^* (G, A, \af)$ is isometrically rich.
\end{cor}

\begin{proof}
Using the condition in Theorem~1.1(ii) of~\cite{BP2},
one sees that quasi-invertible elements in~$A$ are one-sided invertible.
So $A$ is isometrically rich,
and Proposition~\ref{P:IsoRich} applies.
\end{proof}

\begin{qst}\label{Q:ExtR}
Let $A$ be an extremally rich unital  \ca\  (not necessarily prime),
let $G$ be a finite group,
and let $\af \colon G \to \Aut (A)$ have the Rokhlin property.
Does it follow that $C^* (G, A, \af)$ is extremally rich?
\end{qst}

The methods of this paper appear to break down,
because local approximation by extremally rich \ca s
does not necessarily imply extremal richness.
In fact, in Example~5.3 of~\cite{BP2}
there is a direct system of extremally rich \ca s,
with injective unital maps,
such that the direct limit is not extremally rich.

The following result was obtained, using different methods
(and for compact groups),
in Corollary~3.4(2) of~\cite{HW}.
We include the proof here as an example of our methods.

\begin{prp}\label{P:RokhAppDv}
Let $A$ be a separable unital C*-algebra
which is approximately divisible in the sense of~\cite{BKR},
let $G$ be a finite group,
and let $\af \colon G \to \Aut (A)$ have the Rokhlin property.
Then $C^* (G, A, \af)$ is approximately divisible.
\end{prp}

\begin{proof}
Take ${\mathcal{C}}$ to be the collection of all separable
unital approximately divisible \ca s.
This class
has approximate quotients (Definition~\ref{D:Flex}),
because in fact every quotient
of an approximately divisible \ca\  %
is easily seen to be approximately divisible.
So $C^* (G, A, \af)$ is a unital local ${\mathcal{C}}$-algebra
by Theorem~\ref{T:CrPrdRk}.

It is immediate that ${\mathcal{C}}$ is closed under isomorphism.
It is easy to check that ${\mathcal{C}}$ is closed under direct sums
and tensoring with $M_n,$
and ${\mathcal{C}}$ is closed under passing to corners by
Corollary~2.9 of~\cite{BKR}.
So ${\mathcal{C}}$ is finitely saturated.

To show that $C^* (G, A, \af)$ is approximately divisible,
let $F \S C^* (G, A, \af)$ be finite and let $\ep > 0.$
Choose a unital approximately divisible C*-algebra $B$
and a unital homomorphism $\ph \colon B \to A$ such that
$\dist (a, \, \ph (B)) < \frac{1}{3} \ep$ for all $a \in F.$
Choose a finite set $S \S B$ such that
$\dist (a, \, \ph (S)) < \frac{1}{3} \ep$ for all $a \in F.$
Passing to $B / \ker (\ph),$ we may assume $\ph$ is injective.
Choose a unital completely noncommutative finite dimensional \ca\  %
$D \S B$ such that $\| [ b, x] \| < \frac{1}{3} \ep$
for all $x \in S$ and all $b \in D$ with $\| b \| \leq 1.$
Then $\ph (D) \S A$
is a unital completely noncommutative finite dimensional \ca\  %
such that $\| [ b, x] \| < \ep$
for all $x \in F$ and all $b \in \ph (D)$ with $\| b \| \leq 1.$
Thus $C^* (G, A, \af)$ is approximately divisible.
\end{proof}

We next consider the property that the order on \pj s
over~$A$ be determined by traces.

\begin{ntn}\label{N:TraceNtn}
Let $A$ be a unital \ca.
We denote by $T (A)$ the set of all tracial states on $A,$
equipped with the weak* topology.
For any element of $T (A),$
we use the same letter for its standard extension to $M_n (A)$
for arbitrary $n,$
and to $\Mi (A) = \bigcup_{n = 1}^{\infty} M_n (A)$ (no closure).
\end{ntn}

\begin{dfn}\label{D:OrdDetD}
Let $A$ be a simple unital \ca.
We say that the {\emph{order on \pj s over $A$ is determined by traces}}
if whenever $p, q \in \Mi (A)$ are \pj s such that
$\ta (p) < \ta (q)$ for all $\ta \in T (A),$
then $p \precsim q.$
\end{dfn}

This is Blackadar's Second Fundamental Comparability Question
for $\Mi (A).$
See 1.3.1 in~\cite{Bl4}.

\begin{prp}\label{P:OrdDetT}
Let $A$ be a simple separable unital C*-algebra
such that the order on \pj s over $A$ is determined by traces
in the sense of Definition~\ref{D:OrdDetD}.
Let $G$ be a finite group,
and let $\af \colon G \to \Aut (A)$ have the Rokhlin property.
Then the order on \pj s over $C^* (G, A, \af)$ is determined by traces.
\end{prp}

\begin{proof}
The crossed product is simple,
by Remark~1.4 and Corollary~1.6 of~\cite{PhtRp1a}.
So Definition~\ref{D:OrdDetD} applies.

Note that
$M_n (C^* (G, A, \af))
  \cong C^* (G, \, M_n \otimes A, \, \id_{M_n} \otimes \af),$
and $\id_{M_n} \otimes \af$ also has the Rokhlin property.
It therefore suffices to verify the condition
of Definition~\ref{D:OrdDetD} for \pj s in $C^* (G, A, \af).$

Take ${\mathcal{C}} = \{ A \}.$
The finite saturation of ${\mathcal{C}}$ consists of finite
direct sums of unital \ca s which are stably isomorphic to~$A.$
All quotients of such algebras are again algebras of the same form.
So $C^* (G, A, \af)$ is a unital local ${\mathcal{C}}$-algebra
by Theorem~\ref{T:CrPrdRk}.

We finish the proof via an adaptation of a standard argument.
Let $p, q \in C^* (G, A, \af),$
and suppose it is not the case that $p \precsim q.$
We need to find a tracial state on $C^* (G, A, \af)$ such that
$\ta (p) \geq \ta (q).$
Let $\{ x_1, x_2, \ldots \}$ be a countable dense subset
of $C^* (G, A, \af).$
For each $n \geq 1,$
choose a \ca\  $B_n$ in the finite saturation of ${\mathcal{C}}$
and a unital \hm\  $\ph_n \colon B_n \to A$
such that $\dist (x, \, \ph_n (B_n) ) < 2^{-n}$
for all $x \in \{ p, q, x_1, x_2, \ldots, x_n \}.$
{}From the description of the finite saturation of ${\mathcal{C}}$
given above,
we may assume that each $\ph_n$ is injective.
We therefore treat $B_n$ as a unital subalgebra of $C^* (G, A, \af),$
and delete $\ph_n$ from the notation.
We further write
\[
B_n = B_{1, n} \oplus B_{2, n} \oplus \cdots \oplus B_{k (n), \, n}
\]
for unital \ca s $B_{j, n}$ which are stably isomorphic to~$A.$
Each $B_{j, n}$ is a corner of a matrix algebra over~$A,$
so the order on \pj s over $B_{j, n}$ is determined by traces.

Since $\dist (p, B_n), \, \dist (q, B_n) < 2^{-n} \leq \frac{1}{2},$
there exist \pj s $e_n, f_n \in B_n$ such that
$e_n \sim p$ and $f_n \sim q.$
Moreover, we may arrange that $\limi{n} e_n = p$ and $\limi{n} f_n = q.$
Write
\[
e_n = (e_{1, n}, e_{2, n}, \ldots, e_{k (n), \, n})
\andeqn
f_n = (f_{1, n}, f_{2, n}, \ldots, f_{k (n), \, n})
\]
with $e_{j, n}, f_{j, n} \in B_{j, n}.$
It is not the case that $e_n \precsim f_n,$
so there is $j (n)$ such that it is not the case that
$e_{j (n), \, n} \precsim f_{j (n), \, n}.$
Therefore there is a tracial state $\ta_n$ on $B_{j (n), \, n}$
such that $\ta_n (e_{j (n), \, n}) \geq \ta_n (f_{j (n), \, n}).$
Regard $\ta_n$ as a tracial state on $B_n$ by taking $\ta_n = 0$
on the other summands of~$B_n.$
By the Hahn-Banach Theorem,
there is a state $\om_n$ on $C^* (G, A, \af)$
whose restriction to $B_n$ is~$\tau_n.$

Let $\ta$ be a state which is a weak* limit point of the states~$\om_n.$

We claim that $\ta$ is tracial.
Let $a, b \in C^* (G, A, \af),$ and let $\ep > 0.$
We may assume $\ep < 1.$
Set $M = 1 + \max (\| a \|, \| b \|).$
Choose $n_0$ such that $2^{- n_0} < \tfrac{1}{6} M^{-1} \ep$
and so large that there are $k, l \in \{ 1, 2, \ldots, n_0 \}$
such that
\[
\| x_k - a \| < \tfrac{1}{6} M^{-1} \ep
\andeqn
\| x_l - b \| < \tfrac{1}{6} M^{-1} \ep.
\]
Note that $\| x_k \|, \| x_l \| < M.$
Therefore
\[
\| a b - x_k x_l \|
  \leq \| a - x_k \| \cdot \| b \| + \| x_k \| \cdot \| b - x_l \|
  < \big( \tfrac{1}{6} M^{-1} \ep \big) M
        + M \big( \tfrac{1}{6} M^{-1} \ep \big)
  = \tfrac{1}{3} \ep.
\]
For all $n \geq n_0,$ we have $\om_n (x_k x_l) = \om_n (x_l x_k),$
so
\[
| \om_n (a b) - \om_n (b a) |
  \leq 2 \| a b - x_k x_l \|
  < \tfrac{2}{3} \ep.
\]
Choose $n \geq n_0$ such that
\[
| \ta (a b) - \om_n (a b) | < \tfrac{1}{6} \ep
\andeqn
| \ta (b a) - \om_n (b a) | < \tfrac{1}{6} \ep.
\]
Then
\[
| \ta (a b) - \ta (b a) |
 < \tfrac{1}{3} \ep + | \om_n (a b) - \om_n (b a) |
 < \tfrac{1}{3} \ep + \tfrac{2}{3} \ep
 = \ep.
\]
Since $\ep > 0$ is arbitrary, the claim follows.

We complete the proof by showing that $\ta (p) \geq \ta (q).$
Let $\ep > 0.$
Choose $n_0$ such that $n \geq n_0$ implies
\[
\| e_n - p \| < \tfrac{1}{4} \ep
\andeqn
\| f_n - q \| < \tfrac{1}{4} \ep.
\]
Choose $n \geq n_0$ such that
\[
| \ta (p) - \om_n (p) | < \tfrac{1}{4} \ep
\andeqn
| \ta (q) - \om_n (q) | < \tfrac{1}{4} \ep.
\]
Then
\[
| \ta (p) - \om_n (e_n) |
  \leq | \ta (p) - \om_n (p) | + \| p - e_n \|
  < \tfrac{1}{4} \ep + \tfrac{1}{4} \ep
  = \tfrac{1}{2} \ep,
\]
and similarly $| \ta (q) - \om_n (f_n) | < \tfrac{1}{2} \ep.$
Since $\om_n (e_n) \geq \om_n (f_n),$
we conclude that $\ta (p) > \ta (q) - \ep.$
Since $\ep > 0$ is arbitrary, the proof is complete.
\end{proof}

The following result does not use the local approximation
result of Theorem~\ref{T:CrPrdRk}.
However, we have not found it in the literature,
and it is included here because the statement
seems related to other results here.
The assumption that the action have the Rokhlin property is
stronger than necessary, at least when the \ca\  is simple:
the tracial Rokhlin property suffices.
See Proposition~5.7 of~\cite{ELPW}.

\begin{prp}\label{T:RkAndTr}
Let $A$ be a separable unital C*-algebra.
Let $G$ be a finite group,
and let $\af \colon G \to \Aut (A)$ have the Rokhlin property.
Then the restriction map defines a bijection
from $T (C^* (G, A, \af))$ (see Notation~\ref{N:TraceNtn})
to the set $T (A)^G$ of $G$-invariant tracial states on~$A.$
\end{prp}

\begin{proof}
For $g \in G,$ let $u_g \in C^* (G, A, \af)$ be the standard unitary
in the crossed product.
Let $E \colon C^* (G, A, \af) \to A$ be the standard conditional
expectation,
given by $E \left( \sum_{g \in G} a_g u_g \right) = a_1$
when $a_g \in A$ for $g \in G.$

We claim that the map $\ta \mapsto \ta \circ E$
is an inverse of the restriction map.
First, if $\ta \in T (A)^G,$
one easily checks that $\ta \circ E$ is a tracial state
on $C^* (G, A, \af),$
and it is immediate that its restriction to $A$ is~$\ta.$

So suppose $\ta \in T (C^* (G, A, \af)).$
Proving that $\ta = (\ta |_A) \circ E$ is equivalent to proving
that $\ta (a u_g) = 0$ for $g \in G \setminus \{ 1 \}$
and $a \in A.$
Let $\ep > 0.$
Choose \pj s $e_h \in A$ for $h \in G$ according to the Rokhlin property
(Definition~\ref{SRPDfn}),
with
\[
\dt = \frac{\ep}{(1 + \| a \|) \card (G)}
\]
in place of $\ep$ and with $\{ a \}$ in place of~$F.$
For any $h \in G,$
we then have
$e_h u_g e_h = e_h (u_g e_h u_g^* - e_{g h}) u_g + e_h e_{g h} u_g.$
Since $g \neq 1,$ we have $e_h e_{g h} = 0,$
whence $\| e_h u_g e_h \| < \dt.$
Now
\begin{align*}
| \ta (a u_g) |
& \leq \sum_{h \in G} | \ta ( a u_g e_h^2 ) |
  = \sum_{h \in G} | \ta (e_h a u_g e_h) |
  \leq \sum_{h \in G}
     \big( \| e_h a - a e_h \| + | \ta ( a e_h u_g e_h ) | \big)
        \\
& \leq \sum_{h \in G}
     \big( \| e_h a - a e_h \| + \| a \| \cdot \| e_h u_g e_h \| \big)
  < \card (G) (1 + \| a \|) \dt
  = \ep.
\end{align*}
Since $\ep > 0$ is arbitrary, we conclude that $\ta (a u_g) = 0.$
\end{proof}

\begin{cor}\label{C:RokhUniq}
Let $A$ be a separable unital C*-algebra with a unique tracial state.
Let $G$ be a finite group,
and let $\af \colon G \to \Aut (A)$ have the Rokhlin property.
Then $C^* (G, A, \af)$ has a unique tracial state.
\end{cor}

\begin{proof}
Apply Proposition~\ref{T:RkAndTr},
noting that the tracial state on $A$ is necessarily $G$-invariant.
\end{proof}

\end{document}